\documentclass[11pt]{article}
\usepackage{amsmath,amssymb,latexsym}

\newcommand{\beql}[1]{\begin{equation}\label{#1}}
\newcommand{\eeq}{\end{equation}}

\begin{document}
\begin{center}
{\large\bf Eight Hateful Sequences} \\
\vspace*{+.1in}
{\em N. J. A. Sloane} \smallskip \\
Algorithms and Optimization Department  \\
AT\&T Shannon Lab \\
Florham Park, NJ 07932--0971 \medskip \\
Email address: {\tt njas@research.att.com}
\bigskip \\
March 9, 2008; revised May 13, 2008 \bigskip
\end{center}

\setlength{\baselineskip}{1.0\baselineskip}

%\paragraph{ }
\noindent
Dear Martin Gardner: \\
In your July 1974 {\em Scientific American} column you mentioned
the {\em Handbook of Integer Sequences},
which then contained 2372 sequences.
Today the {\em On-Line Encyclopedia of Integer Sequences}
(the {\em OEIS}) \cite{OEIS} contains 140000 sequences.
Here are eight of them, suggested by the theme
of the Eighth Gathering: they are all infinite,
and all \, 'ateful in one way or another. I hope
you like 'em!
Each one is connected with an interesting unsolved problem.

Since this is a $15$-minute talk, I can't give 
many details---see
the entries in the OEIS for more information, and 
for links to related sequences.

\paragraph{1. Hateful or Beastly Numbers}
The most hateful sequence of all!
These are the numbers that contain the string 666 in their
decimal expansion:
$$
666,1666,2666,3666,4666,5666,6660,6661,6662,6663,6664,6665,6666,6667, \ldots \,
$$
(A051003). This sequence is doubly hateful, because it is
based on superstition and because it depends on the fact that
we write numbers in base $10$.

It has been said that if the number of the Beast is $666$, its
fax number must be $667$; fortunately
this has not yet been made the basis for any sequence.
(Added later: unfortunately this is no longer true---see A138563!)

Base-dependent sequences are not encouraged in 
the OEIS, but nevertheless many are present because 
someone has found them interesting, or they have
appeared on a quiz or a web-site, etc.  
One has to admit that some of these sequences 
are very appealing.

For example: start with $n$; if it is a palindrome, stop;
otherwise add $n$ to itself with the digits reversed; repeat until 
you reach a palindrome; or set the value to $-1$ if you
never reach a palindrome. Starting at $19$, we have
$$
19 \rightarrow 19 + 91 = 110 \rightarrow 110 + 011 = 121 \,,
$$
which is a palindrome, so we stop, and the $19$th term
of the sequence is $121$.

The sequence (A033865) begins:
$$
0,1,2,3,4,5,6,7,8,9,11,11,33,44,55,66,77,88,99,121,22,33,22,55,66,77, \ldots \,.
$$
The first unsolved case occurs when we start with $196$. Sequence A006960 
gives its trajectory, which begins
$$
196,887,1675,7436,13783,52514,94039,187088,1067869,10755470,18211171, \ldots \,.
$$
It {\em appears} that it never reaches a palindrome,
but it would be nice to have a proof!
This is hateful because it is one of those problems that seem
too difficult for twenty-first century mathematics to solve.

\paragraph{2.  \'{E}ric Angelini's ``1995'' puzzle }
The following puzzle was invented by  \'{E}ric Angelini
in September 2007: find the
rule that generates the sequence
$$
\mbox{o~n~e~~ n~i~n~e~~ n~i~n~e~~ f~i~v~e~~
f~i~v~e~~ n~i~n~e~~ n~i~n~e~~ f~i~v~e~~
} \cdots \, .
$$
The answer is that if each letter is replaced by its
rank in the Roman alphabet, then the absolute values
of the differences between successive numbers
produce the same sequence:
$$
1,9,9,5,5,9,9,5,5,9,1,3,13,17,1,3,13,17,9,5,5,9,9,5,5,9,1,3,13,17,\ldots
$$
(A131744).
In \cite{1995}, David Applegate and I analyzed
this sequence, and showed among
other things that only $19$ numbers occur
($16$, $19$, $20$ and $22$--$26$ never appear).
We determined the relative frequencies
of these $19$ numbers: $9$ occurs the
most often, with density $0.173\ldots$.

In English this sequence is unique: ``one'' is the only number with
the property that it is equal to 
the absolute value of the difference in rank 
between the first two letters of its name (``o'' is the
fifteenth letter of the alphabet, ``n''
is the fourteenth, and $1 = |14-15|$).
We also discuss versions in other languages.
In French one can begin with either $4$ (see A131745)
or $9$ (A131746), in German with either
$9$ (A133816) or $15$ (A133817), in Italian with 4 (A130316), 
in Russian with either $1$ (A131286) or $2$ (A131287),
and so on.

\paragraph{3. Powertrains }
What is the next term in the following sequence?
$$
679 \rightarrow 378 \rightarrow 168 \rightarrow 48 \rightarrow 32 \rightarrow ?
$$
Answer: 6. The reason is that each term is the product of
the digits of the previous term. Eventually every number $n$
reaches a single-digit number (these are
the only fixed-points), and 
the number of steps for this to
happen is called the {\em persistence} of $n$.
$679$ has persistence $5$, and is in fact the
smallest number with persistence $5$. The smallest numbers
with persistence $n = 1, 2, \ldots, 11$ are given
in sequence A003001:
$$
10,25,39,77,679,6788,68889,2677889,26888999,3778888999,
277777788888899.
$$
This sequence was the subject of an article I wrote in 1973 
\cite{me33}, which Martin Gardner may remember!
I conjectured that the sequence is finite,
and even today no number of persistence greater
than $11$ has been found.

In December 2007, John Conway proposed  some variations
of ``persistence'', one of which I will
discuss here (see \cite{PT} for further
information). 
If $n$ has decimal expansion $abcd\ldots$, the
{\em powertrain} of $n$ is the number
$a^b c^d \ldots$,
which ends in an exponent or a base
according as the number of digits in $n$
is even or odd. We take $0^0 = 1$, and
define the powertrain of $0$ to be $0$.

The OEIS now contains numerous sequences
related to powertrains: see A133500 and
the sequences cross-referenced there.
For example, the following numbers
are fixed under the powertrain map:
$$
0,1,2,3,4,5,6,7,8,9,2592,24547284284866560000000000
$$
(A135385), and we conjecture that there are no others.
Certainly there are no other fixed points below $10^{100}$.

\paragraph{4. Alekseyev's ``$123$'' sequence}
This is a question about strings, proposed by Max Alekseyev
when he was a graduate student in the Computer Science
Department of the University of
California at San Diego (personal communication).
If we start with the string $12$, and repeatedly 
duplicate any substring in place, 
the strings we obtain are:
$$
\begin{array}{ccccccccccccccc}
~ & ~ & ~ & ~ & ~ & ~ & ~ &       12 & ~ & ~ & ~ & ~ & ~ & ~ & ~ \\
~ & ~ & ~ & ~ & ~ & ~ &    112 &  ~ &  122 & ~ & ~ & ~ & ~ & ~ & ~ \\
~ & ~ & ~ & ~ & 1112 & ~ & 1122 & ~ & 1212 & ~ & 1222 & ~ & ~ & ~ & ~ \\
11112 & ~ & 11122 & ~ & 11212 & ~ & 11222 & ~ & 12112 & ~ & 12122 & ~ & 12212 & ~ & 12222 \\
\multicolumn{15}{c}{\cdots}
\end{array}
$$
(A130838). These strings must start with $1$ and end with $2$,
but are otherwise arbitrary. So the number
of such strings of length $n$ is $2^{n-2}$ for $n \ge 2$.

But what if we start with the string $123$?
Now the strings we obtain are
$$
\begin{array}{ccccccccccccccc}
~ & ~ & ~ & ~ & ~ & ~ & ~ &      123 & ~ & ~ & ~ & ~ & ~ & ~ & ~ \\
~ & ~ & ~ & ~ & ~ & 1123 &    ~ &  1223 &  ~ & 1233 & ~ & ~ & ~ & ~ & ~ \\
11123 & ~ & 11223 & ~ & 11233 & ~ & 12223 & ~ & 12233 & ~ & 12333 & ~ & 12123 & ~ & 12323 \\
\multicolumn{15}{c}{\cdots}
\end{array}
$$
%$$
%123,
%$$
%$$
%1123,~1223,~1233,
%$$
%$$
%11123,~11223,~11233,~12223,~12233,~12333,~12123,~12323,
%$$
%$$
%\ldots \,
%$$
(A135475), and the number of such strings of length $n \ge 3$ is
$$
1,3,8,21,54,138,355,924,2432,6461,17301,46657,126656,345972,950611, \ldots \,
$$
(A135473). The question is, what is the $n$-th term in the
latter sequence?
This is hateful because one feels that if this sequence
was only looked at in the right way, there would be a simple recurrence
or generating function.

\paragraph{5. The Curling Number Conjecture} 
Let $S = S_1 S_2 S_3 \cdots S_n$ be a finite string
(the symbols can be anything you like). Write
$S$ in the form $X Y Y \cdots Y = XY^k$,
consisting of a prefix $X$ (which
may be empty), followed
by $k$ (say) copies of a nonempty string $Y$.
In general there will be several ways to do this;
pick one with the greatest value of $k$. Then $k$
is called the {\em curling number} of $S$.

A few years ago Dion Gijswijt proposed the sequence
that is obtained by starting with the string $1$, and 
extending it by continually appending the curling
number of the current string.
The resulting sequence (A090822)
$$
1,1,2,1,1,2,2,2,3,1,1,2,1,1,2,2,2,3,2,1,1,2,1,1,2,2,2,3,1,1,2,1,1, \ldots \,.
$$
was analyzed in \cite{GIJ},
and I talked about it at the Seventh Gathering for Gardner
\cite{G4G7}. It is remarkable because, although it is
unbounded, it grows {\em very} slowly.
For instance, the first $5$ only
appears after about $10^{10^{23}}$ terms.

Some of the proofs in that paper could have
been shortened if we had been able to 
prove a certain conjecture, which remains open to this day.

This is the {\em Curling Number Conjecture},
which states that if one starts with any finite 
string, over any alphabet, and repeatedly
extends it by appending the curling
number of the current string, then eventually
one must reach a $1$.

One way to attack this problem is to start with a string 
that only contains $2$'s and $3$'s, and see how far
we can get before a $1$ appears. For an initial string
of $n$ $2$'s and $3$'s, for $n = 1, 2, \ldots 30$ respectively,
the longest string that can
be obtained before a $1$ appears is:
\begin{align}\label{EqCNC}
1,4,5, & 8,9,14,15,66,68,70,123,124,125,132,133,134,135,136, \nonumber \\
& 138,139, 140,142,143,144,145,146,147,148,149,150,\ldots
\end{align}
(A094004). For example, the best initial string of length $6$ is
$2,2,2,3,2,2$. This produces the string
$$
2,2,2,3,2,2,2,3,2,2,2,3,3,2,1,\ldots 
$$
which extends for $14$ terms before a $1$ appears.
\eqref{EqCNC} is hateful because 
it seems hard to predict the asymptotic behavior. 
Does it continue in a roughly linear
manner for ever, or are there bigger and bigger jumps?
If the Curling Number Conjecture is false, 
the terms could be a lazy eight ($\infty$) from some point on!

\paragraph{6. Leroy Quet's prime-generating recurrence}
Leroy Quet has contributed many original
sequences to the OEIS. Here is one which
so far has resisted attempts to
analyze it. Let $a(0)=2$, and
for $n \ge 1$ define $a(n)$ to
be the smallest prime $p$, not already in
the sequence, such that $n$ divides $a(n-1)+p$.
The sequence begins:
$$
2,3,5,7,13,17,19,23,41,31,29,37,11,67,59,61,83,53,73,79,101,109,89,233,\ldots
$$
(A134204). Is it infinite? To turn the question around, does
$a(n-1)$ ever divide $n$?
If it {\em is} infinite, is it a permutation of the primes?
(This is somewhat reminiscent of the EKG sequence A064413 \cite{EKG}.)
This is hateful because we do not know the answers!

David Applegate has checked that the sequence
exists for at least $450 \cdot 10^6$ terms.
Note that $a(n)$ can be less than $n$. This
happens for the following values of $n$ (A133242):
$$
12,201,379,474,588,868,932,1604,1942,2006,3084,4800,7800,\ldots .
$$

\paragraph{7. The $n$-point traveling salesman problem} 
My colleagues David Johnson and David Applegate have
been re-examining the old question of
the expected length of a traveling salesman tour
through $n$ random points in the unit square
(cf. \cite{TSP},  \cite{BHH}).
To reduce the effects of the boundary, they identify the
edges of the square so as to obtain a flat torus,
and ask for the expected length, $L(n)$,
of the optimal tour through $n$ random points \cite{DSJ08}.
It is convenient to express $L(n)$ in units of eels,
where an {\em eel} is the {\em expected Euclidean
length} of the line joining a random point in the unit
square to the center, a quantity which
is well-known \cite[p.~479]{Finch} to be
$$
\frac{\sqrt{2} + \log(1+\sqrt{2})}{6} ~=~
0.382597858232\ldots  
$$  
(although our name for it is new!).
Then $L(1) = 0$, $L(2) = 2$ eels,
and it is not difficult to show that $L(3) = 3$ eels.
David Applegate has made Monte Carlo estimates
of $L(n)$ for $n \le 50$ based on 
finding optimal tours through a million 
sets of random points. His estimates for $L(2)$ through $L(10)$,
expressed in eels, are:
$$
1.99983,
2.99929,
3.60972,
4.08928,
4.5075,
4.88863,
5.24065,
5.5712,
5.8825,
6.17719.
$$
Observe that the first two terms are a close match for
the true values.
It would be nice to know the exact value of $L(4)$.
The sequence $\{L(n)\}$ is hateful because for $n \ge 4$ it
may be irrational, even when expressed in terms of eels, and so will
be difficult to include in the OEIS.

\paragraph{8. Lagarias's Riemann Hypothesis sequence}
Every inequality in number theory is potentially the source
of a number sequence (if $f(n) \ge g(n)$, set 
$a(n) = \lfloor f(n)-g(n) \rfloor $).
Here is one of the most remarkable examples.
Consider the sequence defined by
$$
a(n) = \lfloor H(n) ~+~ \exp(H(n)) \, \log(H(n))\rfloor ~-~ \sigma(n) \,,
$$
where $H(n) = \sum_{k=1}^{n} \frac{1}{k}$ is the $n$th
harmonic number (see A001008/A002805)
and $\sigma(n)$ is the sum of the divisors of $n$
(A000203).
This begins
$$
0,0,1,0,4,0,7,2,7,5,13,0,17,9,12,8,23,5,27,8,21,20,34,1,33,25, \ldots \,
$$
(A057641).
Jeff Lagarias \cite{Lagarias2002}, extending earlier work of G. Robin
\cite{Robin}, has shown that proving that 
$a(n) \ge 0$ for all $n$ is equivalent to 
proving the Riemann hypothesis!
Hateful because {\em hard}.

\paragraph{ Acknowledgment}
I would like to thank my colleague David Applegate
for being always ready to help analyze a new sequence,
as well as for stepping in on several occasions to keep the OEIS web site up
and running when it was in difficulties.

\end{document}